\documentclass{article}
\usepackage{amsthm}
\usepackage{amssymb}
\usepackage{amsxtra}
 
\newtheorem{thm}{Theorem}[section]
\newtheorem{lemma}{Lemma}[section]
\newtheorem{cor}{Corollary}[section]

\title{More precise Pair Correlation Conjecture}
\author{Tsz Ho Chan}

\begin{document}
\maketitle
\begin{abstract}
In this paper, we derive a more precise version of the Strong Pair Correlation
Conjecture on the zeros of the Riemann zeta function under Riemann Hypothesis
and Twin Prime Conjecture.
\end{abstract}
\section{Introduction}
In the early 1970s, H. Montgomery studied the distribution of the difference
$\gamma - {\gamma}'$ between the imaginary parts of the non-trivial zeros of the
Riemann zeta function. Let
\begin{equation}
\label{basic}
F(x,T) = \mathop{\sum_{0 \leq \gamma \leq T}}_{0 \leq {\gamma}' \leq T}
         x^{i(\gamma - {\gamma}')} w(\gamma - {\gamma}')
\mbox{ and }
w(u) = {4 \over 4 + u^2}.
\end{equation}
Assuming Riemann Hypothesis, he proved in [\ref{M}] that, as $T \rightarrow
\infty$,
$$F(x,T) \sim {T \over 2 \pi} \log{x} + {T \over 2 \pi x^2} (\log{T})^2$$
for $1 \leq x \leq T$ (actually he only proved for $1 \leq x \leq o(T)$ and the
full range was done by Goldston [\ref{G}]). He conjectured that
$$F(x,T) \sim {T \over 2 \pi} \log{T}$$
for $T \leq x$ which is known as the Strong Pair Correlation Conjecture. From
this, one has the (Weak) Pair Correlation Conjecture:
$$\mathop{\sum_{0 < \gamma, \gamma' \leq T}}_{0 < \gamma - \gamma' \leq {2 \pi
\alpha \over \log{T}}} 1 \sim {T\log{T} \over 2 \pi} \int_{0}^{\alpha} 1 -
\Bigl({\sin{\pi u} \over \pi u}\Bigr)^2 du.$$

In [\ref{C1}], the author proved that, assuming Riemann Hypothesis, for any
$\epsilon > 0$,
\begin{equation}
\label{0.1}
\begin{split}
F(x,T) =& {1 \over 2\pi}T\log{x} + {1 \over x^2}\Bigl[{T \over 2\pi}(\log{T 
\over 2\pi})^2 - 2{T \over 2\pi}\log{T \over 2\pi} \Bigr] \\
&+ O(x\log{x}) + O\Bigl({T \over x^{1/2-\epsilon}}\Bigr)
\end{split}
\end{equation}
for $1 \leq x \leq {T \over \log{T}}$. This gives a more precise formula for
$F(x,T)$ in the range $1 \leq x \leq {T \over \log{T}}$. Meanwhile, in
[\ref{C2}], the author derived a more precise Strong Pair Correlation
Conjecture: For every fixed $\epsilon > 0$ and $A \geq 1+\epsilon$,
\begin{equation}
\label{0.2}
F(x,T) = {T \over 2 \pi}\log{T \over 2\pi}- {T \over 2\pi}+ O(T^{1-\epsilon_1})
\end{equation} 
holds uniformly for $T^{1+\epsilon} \leq x \leq T^A$ with some $\epsilon_1 > 0$.
It would be interesting to know how $F(x,T)$ changes from (\ref{0.1}) to
(\ref{0.2}) when $x$ is close to $T$. We have the following
\begin{thm}
\label{theorem1}
Assume Riemann Hypothesis and Twin Prime Conjecture. For any small $\epsilon >
0$ and any integer $M > 2$,
\begin{eqnarray*}
F(x,T) &=& {T \over 2\pi}\log{x} - {4x \over 3\pi}\int_{0}^{T/x} {\sin{v} \over
v} dv + {x^2 \over \pi T} \Bigl(\sum_{h \leq H^{*}} {{\mathfrak S}(h) \over h^2}
\Bigr) \Bigl(1 - \cos{T \over x}\Bigr) \\
&-& {x \over 2\pi}\int_{1}^{\infty} {\sin{Ty \over x} \over y^2} dy + \Bigl({B
\over 2} + {11 \over 12}\Bigr) {x \over \pi} \int_{1}^{\infty} {\sin{Ty \over x}
\over y^4} dy - {4T \over \pi} \int_{1}^{\infty} {f(y) \over y^2} {\sin{{T \over
x}y} \over {T \over x}y} dy \\
&+& {2T \over \pi} \int_{1}^{\infty} {\int_{1}^{y} f(u) du \over y^3} {\sin{{T
\over x}y} \over {T \over x}y} dy + {6T \over \pi} \int_{1}^{\infty} y \int_{y}^
{\infty} {f(u) \over u^4} du {\sin{{T \over x}y} \over {T \over x}y} dy \\
&+& O\Bigl({x^{1+6\epsilon} \over T}\Bigr) + O(x^{{1/2}+7\epsilon}) + O
\Bigl({T \over (\log{T})^{M-2}}\Bigr).
\end{eqnarray*}
for ${T \over (\log{T})^M} \leq x \leq T^{2-\epsilon}$. $B = -C_0 - \log{2\pi}$
and $C_0$ is Euler's constant $0.5772156649..$. $H^{*}$, ${\mathfrak S}(h)$ and
$f(u)$ are defined as in the next section. The implicit constants may depend on
$\epsilon$ and $M$.
\end{thm}
As corollaries of Theorem \ref{theorem1}, we have
\begin{cor}
\label{corollary1}
Assume Riemann Hypothesis and Twin Prime Conjecture. For any integer $M >2$,
$$F(x,T) = {T \over 2\pi}\log{x} + O(x) + O_{M}\Bigl({T \over (\log{T}
)^{M-2}}\Bigr)$$
for ${T \over (\log{T})^M} \leq x \leq T$.
\end{cor}
\begin{cor}
\label{corollary2}
Assume Riemann Hypothesis and Twin Prime Conjecture. For any small $\epsilon >
0$ and any integer $M > 2$,
$$F(x,T) = {T \over 2\pi} \log{T \over 2\pi} - {T \over 2 \pi} + O_{\epsilon}
\Bigl(T \bigl({T \over x}\bigr)^{{1/2}-\epsilon}\Bigr) + O_{\epsilon, M}
\Bigl({T \over (\log{T})^{M-2}} \Bigr)$$
for $T \leq x \leq T^{2-29\epsilon}$.
\end{cor}
\section{Preparations}
We mentioned Twin Prime Conjecture in the previous section. The form needed is
the following: For any $\epsilon > 0$,
$$\sum_{n=1}^{N} \Lambda(n) \Lambda(n+d) = {\mathfrak S}(d) N + O(N^{1/2+
\epsilon})$$
uniformly in $|d| \leq N$. $\Lambda(n)$ is the von Mangoldt lambda function.
${\mathfrak S}(d) = 2\prod_{p>2}\bigl(1-{1 \over (p-1)^2}\bigr) \prod_{p|d, p>2}
{p-1 \over p-2}$ if $d$ is even, and ${\mathfrak S}(d) = 0$ if $d$ is odd. We
also need a lemma concerning ${\mathfrak S}(d)$.
\begin{lemma}
\label{lemma2.1}
For any $\epsilon > 0$,
$$\sum_{k=1}^{h} (h-k) {\mathfrak S}(k) = {1 \over 2} h^2 - {1 \over 2} h\log{h}
+ Ah + O(h^{{1/2}+\epsilon})$$
where $A = {1 \over 2} (1 - C_0 -\log{2\pi})$ and $C_0$ is Euler's constant.
\end{lemma}

Proof: This is a theorem in Montgomery and Soundararajan [\ref{MS}].

Borrowing from [\ref{GGOS}],
$$S_{\alpha}(y) := \sum_{h \leq y} {\mathfrak S}(h) h^{\alpha} - {y^{\alpha +1}
\over \alpha +1} \mbox{ for } \alpha \geq 0,$$
and
$$T_{\alpha}(y) := \sum_{h>y} {{\mathfrak S}(h) \over h^{\alpha}} \mbox{ for }
\alpha>1.$$
Then from [\ref{FG}],
\begin{equation}
\label{1}
S_0(y) = -{1 \over 2}\log{y} + O((\log{y})^{2/3}).
\end{equation}
Suppose $S_0(y) = -{1 \over 2}\log{y} + \epsilon(y)$. By partial
summation,
\begin{equation}
\label{2}
S_{\alpha}(y) = -{y^{\alpha} \over 2\alpha} + \epsilon(y) y^{\alpha} - \alpha
\int_{1}^{y} \epsilon(u) u^{\alpha -1} du + \Bigl({1 \over 2\alpha} + {\alpha
\over \alpha+1}\Bigr),
\end{equation}
and
\begin{equation}
\label{3}
T_{\alpha}(y) = {1 \over (\alpha -1)y^{\alpha -1}} - {\epsilon(y) \over
y^{\alpha}} - {1 \over 2\alpha y^{\alpha}} + \alpha \int_{y}^{\infty}
{\epsilon(u) \over u^{\alpha +1}} du.
\end{equation}
\begin{lemma}
\label{lemma2.2}
For any $\epsilon > 0$,
$$\int_{1}^{y} \epsilon(u) du = {B \over 2}y + O(y^{{1/2}+\epsilon})$$
where $B = -C_0 - \log{2\pi}$ as in the previous section.
\end{lemma}

Proof: By Lemma \ref{lemma2.1},
\begin{eqnarray*}
\int_{1}^{y} \epsilon(u) du &=& \int_{1}^{y} \Bigl(\sum_{h \leq u} {\mathfrak S}
(h) -u+{1 \over 2} \log{u} \Bigr) du \\
&=&\sum_{h \leq y} (y-h) {\mathfrak S}(h) - {1 \over 2}y^2 + {1 \over 2}y\log{y}
- {1 \over 2}y + 1 \\
&=&Ay - {1 \over 2}y + O(y^{{1/2}+\epsilon}) \\
&=& {B \over 2} y + O(y^{{1/2}+\epsilon}).
\end{eqnarray*}

Now, let us define
$$f(y) := \int_{1}^{y} \epsilon(u) - {B \over 2} du.$$
By integration by parts and Lemma \ref{lemma2.2}, one has
\begin{equation}
\label{4}
\int_{1}^{y} \epsilon(u)u du = {B \over 4}y^2 + yf(y) - \int_{1}^{y} f(u) du -
{B \over 4} = {B \over 4}y^2 + O(y^{{3/2}+\epsilon}),
\end{equation}
and
\begin{equation}
\label{5}
\int_{y}^{\infty} {\epsilon(u) \over u^3} du = {B \over 4 y^2} - {f(y) \over
y^3} + 3 \int_{y}^{\infty} {f(u) \over u^4} du = {B \over 4 y^2} + O(y^{-{5/2}+
\epsilon}).
\end{equation}

Next, we are going to define a smooth weight $\Psi_U (t)$. Fix a small positive
real number $\epsilon$ and let $K$ be a large integer depending on $\epsilon$.
Let $M$ be an integer greater than $2$ and $U = (\log{T})^M$. We want $\Psi_U
(t)$ to have support in $[-1/U, 1+1/U]$, $0 \leq \Psi_U (t) \leq 1$, $\Psi_U (t)
= 1$ for $1/U \leq t \leq 1-1/U$, and $\Psi_{U}^{(j)} (t) \ll U^j$ for $j = 1,2,
...,K$.

Let $\Delta = 1 / (2^K U)$. We define a sequence of functions as follow (which
is Vinogradov's construction) :
\begin{eqnarray*}
\chi_0 (t) &=& \left\{ \begin{array}{ll}
1, & \mbox{if $0 \leq t \leq 1$},\\
0, & \mbox{else}.\end{array} \right. \\
\chi_i (t) &=& {1 \over 2 \Delta} \int_{-\Delta}^{\Delta} \chi_{i-1} (t+x) dx
\mbox{ for } i=1,2,...,K+1.
\end{eqnarray*}

Clearly, $0 \leq \chi_i (t) \leq 1$ for $1 \leq i \leq K+1$. One can easily 
check by induction that $\chi_i (t) =1$ for $2^{i-1} \Delta \leq t \leq 1 -
2^{i-1} \Delta$, and $\chi_i (t)=0$ for $t < -2^{i-1} \Delta$ or $t>1+2^{i-1}
\Delta$ for $i=1,2,...,K+1$.
\begin{lemma}
\label{lemma2.3}
$\chi_{i}^{(j)} (t)$ exist and are continuous, and $\chi_{i}^{(j)} (t) \leq 
\Delta^{-j}$ for $0 \leq j \leq i-1$ and $2 \leq i \leq K+1$.
\end{lemma}

Proof: Induction on $i$. First note that $\chi_1(t)$ is continuous because
\begin{eqnarray*}
|\chi_1(t+\delta) - \chi_1(t)| &=& \Big| {1 \over 2\Delta} \int_{-\Delta}^
{\Delta} \chi_0(t+\delta+x) dx - {1 \over 2\Delta} \int_{-\Delta}^{\Delta} 
\chi_0(t+x) dx \Big| \\
&=& \Big| {1 \over 2\Delta} \int_{-\Delta+\delta}^{\Delta+\delta} \chi_0(t+x) 
dx - {1 \over 2\Delta} \int_{-\Delta}^{\Delta} \chi_0(t+x) dx \Big| \\
&=& \Big| {1 \over 2\Delta} \int_{\Delta}^{\Delta+\delta} \chi_0(t+x) dx -
{1 \over 2\Delta} \int_{-\Delta}^{-\Delta+\delta} \chi_0(t+x) dx \Big| \\
&\leq& {\delta \over \Delta}.
\end{eqnarray*}
Similarly,
\begin{eqnarray*}
{\chi_2(t+h) - \chi_2(t) \over h} &=& {1 \over h} \Bigl[ {1 \over 2\Delta} 
\int_{\Delta}^{\Delta+h} \chi_1(t+x) dx - {1 \over 2\Delta} \int_{-\Delta}^
{-\Delta+h} \chi_1(t+x) dx \Bigr] \\
&=& {1 \over 2\Delta} [ \chi_1(t+\Delta+\xi_1) - \chi_1(t-\Delta+\xi_2) ]
\end{eqnarray*}
for some $0 \leq \xi_1, \xi_2 \leq h$ by mean-value theorem.
So $\chi_{2}'(t)$ exists and equals to ${1 \over 2\Delta} [ \chi_1(t+\Delta) - 
\chi_1(t-\Delta) ]$ which is continuous and $\leq {1 \over \Delta}$.
Assume that $\chi_{i}^{(j)} (t)$ are continuous and satisfy 
$\chi_{i}^{(j)} \ll \Delta^{-j}$ for some $2 \leq i \leq K$ and all $0 \leq j
\leq i-1$. Now, for $0 \leq j \leq i-1$, $\chi_{i+1}^{(j)} (t) = {1 \over
2 \Delta} \int_{-\Delta}^{\Delta} \chi_{i}^{(j)} (t+x) dx \leq \Delta^{-j}$ by
induction hypothesis.
For $j=i$,
\begin{eqnarray*}
\chi_{i+1}^{(i)} (t) &=& \lim_{h \rightarrow 0} {\chi_{i+1}^{(i-1)} (t+h) -
\chi_{i+1}^{(i-1)} (t) \over h} \\
&=& \lim_{h \rightarrow 0} {1 \over h} \Bigl[ {1 \over 2\Delta} \int_{\Delta}^
{\Delta+h} \chi_{i}^{(i-1)} (t+x) dx - {1 \over 2\Delta} \int_{-\Delta}^
{-\Delta+h} \chi_{i}^{(i-1)} (t+x) dx \Bigr] \\
&=& {1 \over 2\Delta} [\chi_{i}^{(i-1)} (t+\Delta) - \chi_{i}^{(i-1)} 
(t-\Delta)]
\end{eqnarray*}
which is continuous and $\leq \Delta^{-i}$ by induction hypothesis.
\begin{lemma}
\label{lemma2.4}
$\hat{\chi}_{0} (y) = e^{\pi i y} {\sin{\pi y} \over \pi y}$ and $\hat{\chi}_
{i+1} (y) = \hat{\chi}_{i} (y) {\sin{2 \pi \Delta y} \over 2 \pi \Delta y}$ 
for $0 \leq i \leq K$. Here $\hat{f}(y)$ denotes the inverse Fourier transform
of $f(t)$, $\hat{f}(y) = \int_{-\infty}^{\infty} f(t) e(yt) dt$.
\end{lemma}

Note: We use inverse Fourier transform so that the notation matches with
[\ref{GG}] and [\ref{GGOS}].

Proof: $\hat{\chi}_{0} (y) = \int_{0}^{1} e(yt) dt = {e^{2\pi i y} - 1 \over
2\pi i y} = e^{\pi i y} {\sin{\pi y} \over \pi y}$.
\begin{eqnarray*}
\hat{\chi}_{i+1} (y) &=& \int_{-\infty}^{\infty} \chi_{i+1} (t) e(yt) dt \\
&=& {1 \over 2\Delta} \int_{-\Delta}^{\Delta} \int_{-\infty}^{\infty} \chi_i
(t+x) e(yt) dt dx \\
&=& {1 \over 2\Delta} \int_{-\Delta}^{\Delta} \hat{\chi}_{i} (y) e(-yx) dx \\
&=& {\hat{\chi}_{i} (y) \over 2\Delta} {e(-y\Delta) - e(y\Delta) \over -2\pi i
y} = \hat{\chi}_{i} (y) {\sin{2\pi \Delta y} \over 2 \pi \Delta y}.
\end{eqnarray*}

Now we take $\Psi_{U} (t) = \chi_{K+1} (t)$, then $\Psi_U (t)$ has the 
required properties by the above discussion and Lemma \ref{lemma2.3}. From 
Lemma \ref{lemma2.4}, we know that $\hat{\Psi}_{U} (y) = e^{\pi i y} 
{\sin{\pi y} \over \pi y} ({\sin{2 \pi \Delta y} \over 2 \pi \Delta y})^{K+1}$.
It follows that
\begin{equation}
\label{452}
\begin{split}
Re \hat{\Psi}_{U}(y) &= {\sin{2\pi y} \over 2\pi y} \Bigl({\sin{2\pi \Delta y}
\over 2\pi \Delta y}\Bigr)^{K+1}, \\
\hat{\Psi}_{U} (y) &\ll y^{-K} \mbox{ for } y \gg T^{\epsilon}, \\
\mbox{ and } \hat{\Psi}_{U} (Ty) &\ll T^{-K \epsilon} \mbox{ for } y \gg 
\tau^{-1} \mbox{ where } \tau = T^{1-\epsilon}.
\end{split}
\end{equation}
These are similar to ($18$) and ($19$) in [\ref{GG}]. Also, by Lemma 
\ref{lemma2.3}, it follows from the discussion in [\ref{GG}] that
$$ \hat{\Psi}_{U} (y), \hat{\Psi}'_{U} (y) \ll \mbox{min}\Bigl(1, ({U \over 2 
\pi y})^K \Bigr) $$
which is ($17$) in [\ref{GG}]. Consequently, the results in [\ref{GG}] are 
true with our choice of $\Psi_U (t)$. Moreover, if one follows their arguments 
carefully, one has their Corollaries $1$ $\&$ $2$ (except that the error term 
may need to be modified by a factor of $N^{\epsilon}$) and Theorem $3$ as long
as $\tau = T^{1-\epsilon} \leq x$.

We shall need the following lemmas concerning our weight function $\Psi_U (t)$.
Here we assume $T \Delta \leq x$.
\begin{lemma}
\label{lemma2.5}
For any integer $n \geq 1$,
$$\int_{1}^{\infty} {1 \over y^n} Re \hat{\Psi}_U \Bigl({Ty \over 2\pi x}\Bigr)
dy = {x \over T} \int_{1}^{\infty} {\sin{Ty \over x} \over y^{n+1}} dy +
O\Bigl(K \Delta \log{1 \over \Delta}\Bigr).$$
\end{lemma}

Proof: By a change of variable $v = {Ty \over x}$ and (\ref{452}), the left
hand side
\begin{eqnarray*}
&=& \Bigl({T \over x}\Bigr)^{n-1} \int_{T/x}^{\infty} {1 \over v^n} {\sin{v}
\over v} \Bigl({\sin{\Delta v} \over \Delta v}\Bigr)^{K+1} dv \\
&=& \Bigl({T \over x}\Bigr)^{n-1} \int_{T/x}^{1/\Delta} {\sin{v} \over v^{n+1}}
(1+O(K \Delta^2 v^2)) dv + O\Bigl(\Bigl({T \over x}\Bigr)^{n-1}
\int_{1/\Delta}^{\infty} {1 \over v^{n+1}}dv\Bigr) \\
&=& \Bigl({T \over x}\Bigr)^{n-1} \int_{T/x}^{1/\Delta} {\sin{v} \over v^{n+1}}
dv + O\Bigl(\Bigl({T \over x}\Bigr)^{n-1} K \Delta^2 \int_{T/x}^{1/\Delta} {1
\over v^{n-1}}dv\Bigr) + O\Bigl(\Bigl({T \over x}\Bigr)^{n-1} \Delta^n\Bigr) \\
&=& \Bigl({T \over x}\Bigr)^{n-1} \int_{T/x}^{\infty} {\sin{v} \over v^{n+1}}
dv + O\Bigl(K\Delta \log{1 \over \Delta}\Bigr) \\
&=& {x \over T} \int_{1}^{\infty} {\sin{Ty \over x} \over y^{n+1}} dy +
O\Bigl(K \Delta \log{1 \over \Delta}\Bigr)
\end{eqnarray*}
because $T\Delta \leq x$. Note that the error term comes from the case $n=2$. If
$n \not= 2$, then we can replace the error term by $O(K \Delta)$.
\begin{lemma}
\label{lemma2.6}
If $F(y) \ll y^{-{3/2} + \epsilon}$ for $y \geq 1$, then
$$\int_{1}^{\infty} F(y) Re \hat{\Psi}_U \Bigl({Ty \over 2 \pi x}\Bigr) dy =
\int_{1}^{\infty} F(y) {\sin{{T \over x}y} \over {T \over x}y} dy + O(K \Delta)
.$$
\end{lemma}

Proof: By a change of variables $v={Ty \over x}$ and (\ref{452}), the left hand
side
\begin{eqnarray*}
&=&{x \over T}\int_{T/x}^{\infty} F({x \over T}v) {\sin{v} \over v}
\Bigl({\sin{\Delta v} \over \Delta v}\Bigr)^{K+1} dv \\
&=&{x \over T}\int_{T/x}^{1/\Delta} F({x \over T}v) {\sin{v} \over v} \bigl(1 +
O(K \Delta^2 v^2)\bigr) dv + O\Bigl({x \over T \Delta^{K+1}} \int_{1/\Delta}^
{\infty} {|F({x \over T}v)| \over v^{K+2}} dv \Bigr) \\
&=&{x \over T}\int_{T/x}^{1/\Delta} F({x \over T}v) {\sin{v} \over v} dv +
O\Bigl(K \bigl({T \over x}\bigr)^{{1/2}-\epsilon} \Delta^{{3/2}-\epsilon}\Bigr)
\\
&=&{x \over T}\int_{T/x}^{\infty} F({x \over T}v) {\sin{v} \over v} dv +
O\Bigl({x \over T} \int_{1/\Delta}^{\infty} {|F({x \over T}v)| \over v}dv \Bigr)
+ O(K \Delta) \\
&=&\int_{1}^{\infty} F(y) {\sin{{T \over x}y} \over {T \over x}y} dy + O(K
\Delta).
\end{eqnarray*}
Finally, we need the following
\begin{lemma}
\label{lemma2.7}
Assume Riemann Hypothesis. For any $\epsilon > 0$,
$$\sum_{n \leq x} \Lambda(n)^2 n = {1 \over 2}x^2\log{x} - {1\over 4} x^2 +
O({x^{3/2+\epsilon}})$$
$$\sum_{n >x} {\Lambda(n)^2 \over n^3} = {1 \over 2}{\log{x} \over x^2} + {1
\over 4}{1 \over x^2} + O({1 \over x^{5/2-\epsilon}})$$
where the implicit constants may depend on $\epsilon$.
\end{lemma}

Proof: By partial summation and the form of prime number theorem under Riemann
Hypothesis.
\section{Proof of main results}
Throughout this section, we assume $\tau = T^{1-\epsilon} \leq {T \over 
(\log{T})^M} \leq x$, $U = (\log{T})^M$ for $M>2$, $H^{*}=\tau^{-2} x^{2/(1-
\epsilon)}$, and $\Psi_{U}(t)$ is defined as in the previous section. Keep in
mind the $\epsilon$ and $M$ dependency in the error terms.

\bigskip

Proof of Theorem \ref{theorem1}: Our method is that of Goldston and Gonek
[\ref{GG}]. Let $s = \sigma + it$,
$$A(s) := \sum_{n \leq x} {\Lambda(n) \over n^s} \mbox{ and } A^{*}(s) :=
\sum_{n>x} {\Lambda(n) \over n^s}.$$
Assume Riemann Hypothesis, it follows from Theorem $3.1$ of [\ref{C1}] with
slight modification that
$$F(x,T) =$$
$${1 \over 2\pi} \int_{0}^{T} \Big| {1 \over x} \Bigl( A(-{1 \over 2}
+it) - \int_{1}^{x} u^{1/2-it} du \Bigr) + x \Bigl( A^{*}({3 \over 2}+it) -
\int_{x}^{\infty} u^{-3/2-it} du \Bigr) \Big|^2 dt$$
$$ + O((\log{T})^3).$$
Inserting $\Psi_{U} (t/T)$ into the integral and extending the range of
integration to the whole real line, we can get
\begin{equation}
\label{4.6.1}
F(x,T) = {1 \over 2 \pi x^2} I_1 (x, T) + {x^2 \over 2 \pi} I_2 (x, T) + O
\Bigl({T(\log{T})^2 \over U}\Bigr) + O \Bigl({x^{1+6\epsilon} \over T}\Bigr)
\end{equation}
where
$$I_1(x,T) = \int_{-\infty}^{\infty} \Psi_{U}\Bigl({t \over T}\Bigr) \Big| A(-{1
\over 2}+it) - \int_{1}^{x} u^{1/2-it} du \Big|^2 dt,$$
and
$$I_2(x,T) = \int_{-\infty}^{\infty} \Psi_{U}\Bigl({t \over T}\Bigr) \Big| A^{*}
({3 \over 2}+it) - \int_{x}^{\infty} u^{-3/2-it} du \Big|^2 dt.$$
This is essentially by Lemma $1$ of [\ref{GGOS}] with modification that $V=-{T
\over U}$ and $T-{T \over U}$, and $W = {2T \over U}$. Riemann Hypothesis is
assumed here so that the contribution from the cross term is estimated via
Theorem $3$ of [\ref{GG}].

\bigskip

Now, we assume the Twin Prime Conjecture in the previous section. By Corollary
$1$ of [\ref{GG}] (see also the calculations at the end of [\ref{GG}] and
[\ref{GGOS}]) and Lemma \ref{lemma2.7}, one has,
\begin{eqnarray*}
I_1(x,T) &=& \hat{\Psi}_{U}(0) T \sum_{n \leq x} \Lambda^2 (n) n \\
& &+ 4 \pi \Bigl({T \over 2 \pi}\Bigr)^3 \int_{T/2\pi x}^{\infty} \Bigl(
\sum_{h \leq 2 \pi x v /T} {\mathfrak S} (h) h^2 \Bigr) Re \hat{\Psi}_{U}(v) 
{dv \over v^3} \\
& &- 4\pi \Bigl({T \over 2 \pi}\Bigr)^3 \int_{T/2\pi \tau x}^{\infty} \Bigl(
\int_{0}^{2\pi x v/T} u^2 du \Bigr) Re \hat{\Psi}_{U}(v) {dv \over v^3} \\
& &+ O\Bigl({x^{3+6\epsilon} \over T}\Bigr) + O(x^{{5/2}+7\epsilon}) \\
&=&{1 \over 2}Tx^2 \log{x} - {1 \over 4}Tx^2 \\
& &+4 \pi \Bigl({T \over 2 \pi}\Bigr)^3 \int_{T/2\pi x}^{\infty} \Bigl(
\sum_{h \leq 2 \pi x v /T} {\mathfrak S} (h) h^2 - \int_{0}^{2\pi x v/T} u^2 du
\Bigr) Re \hat{\Psi}_{U}(v) {dv \over v^3} \\
& &- {4\pi \over 3} x^3 \int_{T/2\pi\tau x}^{T/2\pi x} Re \hat{\Psi}_U (v) dv
+ O\Bigl({x^{3+6\epsilon} \over T}\Bigr) + O(x^{{5/2}+7\epsilon}) \\
&=&{1 \over 2}Tx^2 \log{x} - {1 \over 4}Tx^2 \\
& &+4 \pi \Bigl({T \over 2 \pi}\Bigr)^3 \int_{T/2\pi x}^{\infty} \Bigl(
\sum_{h \leq 2 \pi x v /T} {\mathfrak S} (h) h^2 - \int_{0}^{2\pi x v/T} u^2 du
\Bigr) Re \hat{\Psi}_{U}(v) {dv \over v^3} \\
& &-{2 \over 3} x^3 \int_{0}^{T/x} {\sin{v} \over v} dv + O\Bigl({K T x^2 \over
(\log{T})^M}\Bigr) + O\Bigl({x^{3+6\epsilon} \over T}\Bigr)
+ O(x^{{5/2}+7\epsilon}) 
\end{eqnarray*}
because, from (\ref{452}),
\begin{eqnarray*}
\int_{T/2\pi\tau x}^{T/2\pi x} Re \hat{\Psi}_U (v) dv &=& \int_{0}^{T/2\pi x}
{\sin{2\pi v} \over 2\pi v} \Bigl({\sin{2\pi \Delta v} \over 2\pi \Delta
v}\Bigr)^{K+1} dv + O\Bigl({T \over \tau x}\Bigr) \\
&=&{1 \over 2\pi} \int_{0}^{T/x} {\sin{u} \over u} (1 + O(K \Delta^2 u^2)) du
+ O\Bigl({T \over \tau x}\Bigr) \\
&=&{1 \over 2\pi} \int_{0}^{T/x} {\sin{u} \over u} du + O\Bigl({K \Delta^2 T^2
\over x^2}\Bigr) + O\Bigl({T \over \tau x}\Bigr).
\end{eqnarray*}

Similarly, by Corollary $2$ of [\ref{GG}] and Lemma \ref{lemma2.7},
\begin{eqnarray*}
I_2(x,T) &=& \hat{\Psi}_{U}(0) T \sum_{x<n} {\Lambda^2(n) \over n^3} \\
& &+{8 \pi^2 \over T} \int_{0}^{T/2 \pi x} \Bigl(\sum_{1 \leq h \leq H^{*}}
{{\mathfrak S}(h) \over h^2} \Bigr) Re \hat{\Psi}_{U} (v) v dv \\
& &+{8 \pi^2 \over T} \int_{T/2\pi x}^{T H^{*}/2 \pi x} \Bigl(\sum_{2\pi xv/T <
h \leq H^{*}} {{\mathfrak S}(h) \over h^2} \Bigr) Re \hat{\Psi}_{U} (v) v dv \\
& &-{8 \pi^2 \over T} \int_{0}^{T H^{*}/2\pi x} \Bigl(\int_{2\pi xv/T}^{H^{*}}
u^{-2} du\Bigr) Re \hat{\Psi}_{U} (v) v dv \\
& &+O(T^{-1} x^{-1+6\epsilon}) + O(x^{-{3/2} + 6\epsilon}) + O(T^{1-{\epsilon/2}
} x^{-2}) \\
&=&{T \log{x}\over 2x^2} + {1 \over 4} {T \over x^2} \\
& &+ {8 \pi^2 \over T} \int_{0}^{T H^{*}/2\pi x} \Bigl(\sum_{2\pi xv/T <
h \leq H^{*}} {{\mathfrak S}(h) \over h^2} - \int_{2\pi xv/T}^{H^{*}} {du \over
u^2} \Bigr) Re \hat{\Psi}_{U} (v) v dv \\
& & + O\Bigl({x^{-1+6\epsilon} \over T}\Bigr) + O(x^{-{3/2}+ 6\epsilon}).
\end{eqnarray*}
Therefore, by a change of variable $y={2\pi xv \over T}$ and putting back to
(\ref{4.6.1}), we have
\begin{eqnarray*}
F(x,T) &=& {T \over 2\pi}\log{x} + {T \over \pi}\int_{1}^{\infty} \Bigl(
\sum_{h \leq y} {\mathfrak S}(h)h^2 - {y^3 \over 3} \Bigr) Re \hat{\Psi}_{U}
\Bigl({Ty \over 2\pi x}\Bigr) {dy \over y^3} \\
& &+{T \over \pi}\int_{1}^{H^{*}} \Bigl(\sum_{y<h \leq H^{*}} {{\mathfrak S}(h)
\over h^2} - \int_{y}^{H^{*}} {du \over u^2} \Bigr) Re \hat{\Psi}_{U}
\Bigl({Ty \over 2\pi x}\Bigr) y dy \\
& & + {T \over \pi} \int_{0}^{1} \Bigl(\sum_{h \leq H^{*}} {{\mathfrak S}(h)
\over h^2} - \int_{y}^{H^{*}} {du \over u^2} \Bigr) Re \hat{\Psi}_{U}
\Bigl({Ty \over 2\pi x}\Bigr) y dy \\
& &-{x \over 3\pi} \int_{0}^{T/x} {\sin{v} \over v} dv + O\Bigl({KT \over 
(\log{T})^{M-2}}\Bigr) + O\Bigl({x^{1+6\epsilon} \over T}\Bigr) + O(x^{{1/2}+
7\epsilon})
\end{eqnarray*}
\begin{equation}
\label{4.6.main}
= {T \over 2\pi}\log{x} + {T \over \pi}I_1 + {T \over \pi}I_2 -{4x \over 3\pi}
\int_{0}^{T/x} {\sin{v} \over v} dv +
\end{equation}
$${x^2 \over \pi T} \Bigl(\sum_{h \leq H^{*}} {{\mathfrak S}(h) \over h^2}\Bigr)
\Bigl(1 - \cos{T \over x}\Bigr)+ O\Bigl({KT \over (\log{T})^{M-2}}\Bigr) + 
O\Bigl({x^{1+6\epsilon} \over T}\Bigr) + O(x^{{1/2}+7\epsilon}),$$
where $I_1$ and $I_2$ are the first and second integral respectively. This is
because
\begin{eqnarray*}
& &{T \over \pi} \int_{0}^{1} \Bigl(\sum_{h \leq H^{*}} {{\mathfrak S}(h)
\over h^2} - \int_{y}^{H^{*}} {du \over u^2} \Bigr) Re \hat{\Psi}_{U}
\Bigl({Ty \over 2\pi x}\Bigr) y dy \\
&=&{4\pi x^2 \over T} \Bigl(\sum_{h \leq H^{*}} {{\mathfrak S}(h) \over h^2}
\Bigr) \int_{0}^{T/2\pi x} Re \hat{\Psi}_U (v) vdv - {x \over \pi} \int_{0}^
{T/x} {\sin{u} \over u} du + O\Bigl({KT \over (\log{T})^M}\Bigr) \\
&=&{x^2 \over \pi T} \Bigl(\sum_{h \leq H^{*}} {{\mathfrak S}(h) \over h^2}
\Bigr) \int_{0}^{T/x} \sin{u} \bigl(1+ O(K \Delta^2 u^2)\bigr) du \\
& &- {x \over \pi} \int_{0}^{T/x} {\sin{u} \over u} du + O\Bigl({KT \over 
(\log{T})^M}\Bigr) \\
&=&{x^2 \over \pi T} \Bigl(\sum_{h \leq H^{*}} {{\mathfrak S}(h) \over h^2}
\Bigr) \Bigl(1 - \cos{T \over x}\Bigr) - {x \over \pi} \int_{0}^{T/x} {\sin{u}
\over u} du + O\Bigl({KT \over (\log{T})^M}\Bigr)
\end{eqnarray*}
by a similar calculation as before and $T\Delta \leq x$.
With the notation of $S_{\alpha}(y)$ and $T_{\alpha}(y)$,
\begin{eqnarray*}
I_1 &=& \int_{1}^{\infty} S_2(y) Re \hat{\Psi}_{U} \Bigl({Ty \over 2\pi
x}\Bigr) {dy \over y^3} \\
&=&\int_{1}^{\infty} \Bigl[ {-1 \over 4y} + {\epsilon(y) \over y} - {2
\int_{1}^{y} \epsilon(u)u du \over y^3} \Bigr] Re \hat{\Psi}_{U} \Bigl(
{Ty \over 2\pi x}\Bigr) dy \\
& &- 2\int_{1}^{\infty} {f(y) \over y^2} Re \hat{\Psi}_{U} \Bigl({Ty \over 2\pi
x}\Bigr) dy + 2\int_{1}^{\infty} {\int_{1}^{y} f(u) du \over y^3} Re \hat{\Psi}_
{U} \Bigl({Ty \over 2\pi x}\Bigr) dy \\
& &+ \Bigl({B \over 2} + {11 \over 12}\Bigr) {x \over T} \int_{1}^{\infty}
{\sin{Ty \over x} \over y^4} dy + O(K \Delta) \\
&=&\int_{1}^{\infty} \Bigl[ {-1 \over 4y} + {\epsilon(y) \over y} - {2
\int_{1}^{y} \epsilon(u)u du \over y^3} \Bigr] Re \hat{\Psi}_{U} \Bigl(
{Ty \over 2\pi x}\Bigr) dy \\
& &- 2\int_{1}^{\infty} {f(y) \over y^2} {\sin{{T \over x}y} \over
{T \over x}y} dy + 2\int_{1}^{\infty} {\int_{1}^{y} f(u) du \over y^3} {\sin{{T
\over x}y} \over {T \over x}y} dy \\
& &+ \Bigl({B \over 2} + {11 \over 12}\Bigr) {x \over T} \int_{1}^{\infty}
{\sin{Ty \over x} \over y^4} dy + O(K \Delta)
\end{eqnarray*}
by (\ref{2}), (\ref{4}), Lemma \ref{lemma2.5} and Lemma \ref{lemma2.6}. As for
$I_2$, note that by (\ref{1}) and (\ref{3}),
\begin{equation}
\label{4.6.7}
T_2(z) = {1 \over z} + O\Bigl({(\log{z})^{2/3} \over z^2}\Bigr),
\end{equation}
and
\begin{eqnarray*}
\sum_{y<h \leq H^{*}} {{\mathfrak S}(h) \over h^2} - \int_{y}^{H^{*}} {du \over
u^2}
&=&T_2(y)-T_2(H^{*}) -{1 \over y} -{1 \over H^{*}} \\
&=&T_2(y)-{1 \over y} +O\Bigl({(\log{H^{*}})^{2/3} \over (H^{*})^2} \Bigr).
\end{eqnarray*}
Therefore,
\begin{eqnarray*}
I_2 &=& \int_{1}^{H^{*}} \Bigl(T_2(y) - {1 \over y} + O\Bigl({(\log{H^{*}})^
{2/3} \over (H^{*})^2} \Bigr) \Bigr) Re \hat{\Psi}_{U} \Bigl({Ty \over 2\pi x}
\Bigr) y dy \\
&=&\int_{1}^{H^{*}}\Bigl(T_2(y) - {1 \over y}\Bigr) Re \hat{\Psi}_{U} \Bigl({Ty
\over 2\pi x}\Bigr) y dy \\
& &+ O\Bigl({(\log{H^{*}})^{2/3} x^2 \over (H^{*})^2 T^2} \int_{T/2 \pi x}^{T H
^{*} /2 \pi x} |\hat{\Psi}_{U}(v)| vdv \Bigr) \\
&=&\int_{1}^{\infty}\Bigl(T_2(y) - {1 \over y}\Bigr)Re \hat{\Psi}_{U} \Bigl({Ty
\over 2\pi x}\Bigr) y dy + O\Bigl({1 \over T^{\epsilon}}\Bigr)
\end{eqnarray*}
because of (\ref{4.6.7}) and the formula of $\hat{\Psi}_{U}(y)$ in the previous
section that the integral $\int_{H^{*}}^{\infty} \ll {x (\log{H^{*}})^{2/3}
\over T H^{*}} \ll {1 \over T^{\epsilon}}$ by the definition of $H^{*}$ (similar
estimation for the error term). Applying (\ref{3}), (\ref{5}) and Lemma
\ref{lemma2.6},
\begin{eqnarray*}
I_2 &=& \int_{1}^{\infty} \Bigl[{-1 \over 4y} - {\epsilon(y) \over y} + 2y
\int_{y}^{\infty} {\epsilon(u) \over u^3} du \Bigr] Re \hat{\Psi}_{U} \Bigl({Ty
\over 2\pi x}\Bigr) dy + O\Bigl({1 \over T^{\epsilon}}\Bigr) \\
&=&\int_{1}^{\infty} \Bigl[{-1 \over 4y} - {\epsilon(y) \over y} + {B \over 2y}
\Bigr] Re \hat{\Psi}_{U} \Bigl({Ty \over 2\pi x}\Bigr) dy -2\int_{1}^{\infty}
{f(y) \over y^2} Re \hat{\Psi}_{U} \Bigl({Ty \over 2\pi x}\Bigr) dy \\
& & + 6\int_{1}^{\infty} y \int_{y}^{\infty} {f(u) \over u^4} du Re \hat
{\Psi}_{U} \Bigl({Ty \over 2\pi x}\Bigr) dy + O\Bigl({1 \over T^{\epsilon}}
\Bigr) \\
&=&\int_{1}^{\infty} \Bigl[{-1 \over 4y} - {\epsilon(y) \over y} + {B \over 2y}
\Bigr] Re \hat{\Psi}_{U} \Bigl({Ty \over 2\pi x}\Bigr) dy \\
& &-2\int_{1}^{\infty} {f(y) \over y^2} {\sin{{T \over x}y} \over {T \over x}y}
dy + 6\int_{1}^{\infty} y \int_{y}^{\infty} {f(u) \over u^4} du  {\sin{{T \over
x}y} \over {T \over x}y} dy + O(K\Delta).
\end{eqnarray*}
Consequently, with miraculous cancellations, one has
\begin{eqnarray*}
I_1 + I_2 &=& -{x \over 2T} \int_{1}^{\infty} {\sin{Ty \over x} \over y^2} dy +
\Bigl({B \over 2} + {11 \over 12}\Bigr) {x \over T} \int_{1}^{\infty} {\sin{Ty
\over x} \over y^4} dy \\
& &-4\int_{1}^{\infty} {f(y) \over y^2} {\sin{{T \over x}y} \over {T \over x}y}
dy +2\int_{1}^{\infty} {\int_{1}^{y} f(u) du \over y^3} {\sin{{T \over x}y}
\over {T \over x}y} dy \\
& &+6\int_{1}^{\infty} y \int_{y}^{\infty} {f(u) \over u^4} du {\sin{{T \over x}
y} \over {T \over x}y} dy + O(K \Delta)
\end{eqnarray*}
by Lemma \ref{lemma2.5} again. Putting this back to (\ref{4.6.main}), we have
Theorem \ref{theorem1}.

\bigskip

Proof of Corollary \ref{corollary1}: This follows from Theorem \ref{theorem1}
straighforwardly as $x \leq T$ and $f(u) \ll u^{1/2+\epsilon}$ by Lemma
\ref{lemma2.2}. Note that the error term is better than (\ref{0.1}) for $x$ in
the given range.

\bigskip

Before proving Corollary \ref{corollary2}, we need the following lemmas.
\begin{lemma}
\label{4lemma7.1}
$$\int_{1}^{\infty} {\sin{ax} \over x^{2n}} dx = {a^{2n-1} \over (2n-1)!}
\Bigl[\sum_{k=1}^{2n-1} {(2n-k-1)! \over a^{2n-k}} \sin{\bigl(a+(k-1){\pi \over
2}\bigr)} + (-1)^n ci(a) \Bigr]$$
where $ci(x) = -\int_{x}^{\infty} {\cos{t} \over t}dt = C_0 + \log{x} +
\int_{0}^{x} {\cos{t}-1 \over t} dt$ and $C_0$ is Euler's constant.
\end{lemma}

Proof: This is formula $3.761(3)$ on P.430 of [\ref{GR}] which can be proved by
integration by parts inductively.

\begin{lemma}
\label{4lemma7.2}
If $F(y) \ll y^{-{3/2}+\epsilon}$ for $y \geq 1$, then for $T \leq x$,
$$\int_{1}^{\infty} F(y) {\sin{{T \over x}y} \over {T \over x}y} dy =
\int_{1}^{\infty} F(y) dy + O\Bigl(\bigl({T \over x}\bigr)^{{1/2}-\epsilon}
\Bigr).$$
\end{lemma}

Proof: Since $T \leq x$, the left hand side
\begin{eqnarray*}
&=& \int_{1}^{x/T} F(y) \Bigl(1 + O\Bigl(\bigl({T \over x}\bigr)^2 y^2 \Bigr)
\Bigr) dy + O\Bigl(\int_{x/T}^{\infty} {|F(y)| \over {T \over x}y} dy \Bigr) \\
&=& \int_{1}^{x/T} F(y) dy + O\Bigl(\bigl({T \over x}\bigr)^{{1/2}-\epsilon}
\Bigr) \\
&=& \int_{1}^{\infty} F(y) dy + O\Bigl(\int_{x/T}^{\infty} |F(y)| dy \Bigr) + 
O\Bigl(\bigl({T \over x}\bigr)^{{1/2}-\epsilon}\Bigr) \\
&=& \int_{1}^{\infty} F(y) dy + O\Bigl(\bigl({T \over x}\bigr)^{{1/2}-\epsilon}
\Bigr).
\end{eqnarray*}
\begin{lemma}
\label{4lemma7.3}
$$\sum_{h \leq H^{*}} {{\mathfrak S}(h) \over h^2} = {7 \over 4}+ {B \over 2}+
6\int_{1}^{\infty} {f(u) \over u^4} du + O\Bigl({1 \over H^{*}}\Bigr)$$
where $B = -C_0-\log{2\pi}$ and $C_0$ is Euler's constant again.
\end{lemma}

Proof: First, from (\ref{1}),
\begin{eqnarray*}
\sum_{h > H^{*}} {{\mathfrak S}(h) \over h^2} &=& \int_{H^{*}}^{\infty} {1 \over
u^2} d(S_0(u) + u) \\
&\ll& {1 \over H^{*}} + \int_{H^{*}}^{\infty} {\log{u} \over u^3} du \ll {1
\over H^{*}}
\end{eqnarray*}
which accounts for the error term. It remains to see that
\begin{eqnarray*}
\sum_{h=1}^{\infty} {{\mathfrak S}(h) \over h^2} &=& \int_{1}^{\infty} {1 \over
u^2} d(S_0(u)+u) \\
&=&2 \int_{1}^{\infty} {S_0(u)+u \over u^3} du \\
&=&2 \int_{1}^{\infty} {u - {1 \over 2}\log{u} + \epsilon(u) \over u^3} du \\
&=&2 - {1 \over 4} + 2 \int_{1}^{\infty} {\epsilon(u) \over u^3} du \\
&=&{7 \over 4} + {B \over 2} + 2 \int_{1}^{\infty} {\epsilon(u) - {B \over 2}
\over u^3} du \\
&=&{7 \over 4} + {B \over 2} + 2 \int_{1}^{\infty} {1 \over u^3} df(u) \\
&=&{7 \over 4} + {B \over 2} + 6\int_{1}^{\infty} {f(u) \over u^4} du
\end{eqnarray*}
by integration by parts and the definitions of $\epsilon(u)$ and $f(u)$.

\bigskip

Proof of Corollary \ref{corollary2}: First observe that when $x$ is in the
required range, the error terms in Theorem \ref{theorem1} is $O_{\epsilon, M}
\Bigl({T \over (\log{T})^{M-2}}\Bigr)$. Rewrite Theorem \ref{theorem1} as
$$F(x,T)= {T \over 2\pi} \log{x} - T_1 + T_2 - T_3 + T_4 - T_5 + T_6 + T_7 + 
\mbox{error}.$$
Then, by Lemma \ref{4lemma7.1}, Lemma \ref{4lemma7.2} and Lemma \ref{4lemma7.3},
\begin{eqnarray*}
T_1 &=& {4x \over 3\pi}\int_{0}^{T/x} 1+O(v^2) dv = {4T \over 3\pi} + O\Bigl({
T^3 \over x^2}\Bigr), \\
T_2 &=& {x^2 \over \pi T}\Bigl({7 \over 4} + {B \over 2} + 6\int_{1}^{\infty}
{f(u) \over u^4} du + O ({1 \over H^{*}})\Bigr) \Bigl[{1 \over 2}\bigl({T \over
x}\bigr)^2 + O\Bigl(\bigl({T \over x}\bigr)^4\Bigr) \Bigr] \\
&=&{T \over 2\pi} \Bigl({7 \over 4} + {B \over 2} + 6\int_{1}^{\infty}
{f(u) \over u^4} du \Bigr) + O(T^{1-2\epsilon}) + O\Bigl(T \bigl({T \over
x}\bigr)^2 \Bigr), \\
T_3 &=& {T \over 2\pi} \Bigl[{1 \over {T/x}}\sin{T \over x} - ci({T \over x})
\Bigr] \\
&=& -{T \over 2\pi}\log{T \over x} - {C_0 T \over 2\pi} + {T \over 2\pi} +
O\Bigl(T({T \over x})\Bigr), \\
T_4 &=& \Bigl({B \over 2}+{11 \over 12}\Bigr) {x \over 6 \pi} \Bigl({T \over x}
\Bigr)^3 \Bigl[2\bigl({x \over T}\bigr)^3 \sin{T \over x} + \bigl({x \over T}
\bigr)^2 \sin{({T \over x} + {\pi \over 2})} \\
& &+ \bigl({x \over T}\bigr)\sin{({T \over x} + \pi)} + ci({T \over x})
\Bigr] \\
&=& \Bigl({B \over 2}+{11 \over 12}\Bigr) {T \over 2\pi} + O\Bigl(T \bigl({T 
\over x}\bigr)\Bigr),
\end{eqnarray*}
\begin{eqnarray*}
T_5 &=& {4T \over \pi} \int_{1}^{\infty} {f(y) \over y^2} dy + O\Bigl( T
\bigl({T \over x}\bigr)^{{1/2}-\epsilon} \Bigr), \\
T_6 &=& {2T \over \pi} \int_{1}^{\infty} {\int_{1}^{y} f(u) du \over y^3} dy +
O\Bigl( T\bigl({T \over x}\bigr)^{{1/2}-\epsilon} \Bigr) \\
&=& {T \over \pi} \int_{1}^{\infty} {f(y) \over y^2} dy + O\Bigl( T\bigl({T
\over x}\bigr)^{{1/2}-\epsilon} \Bigr), \\
T_7 &=& {6T \over \pi} \int_{1}^{\infty} y \int_{y}^{\infty} {f(u) \over u^4}
du dy + O\Bigl( T\bigl({T \over x}\bigr)^{{1/2}-\epsilon} \Bigr) \\
&=&-{3T \over \pi} \int_{1}^{\infty} {f(u) \over u^4} du + {3T \over \pi} 
\int_{1}^{\infty} {f(y) \over y^2} dy + O\Bigl( T\bigl({T \over x}\bigr)^{{1/2}
-\epsilon} \Bigr).
\end{eqnarray*}
Combining these, we get
\begin{eqnarray*}
F(x,T) &=& {T \over 2\pi}\log{T} + {T \over 2\pi} \Bigl[-{8 \over 3} + {7 \over
4} + {B \over 2} + C_0 -1 + {B \over 2} + {11 \over 12}\Bigr] \\
& & + O\Bigl( T\bigl({T \over x}\bigr)^{{1/2}-\epsilon} \Bigr) + O\Bigl({T \over
(\log{T})^{M-2}} \Bigr) \\
&=& {T \over 2\pi}\log{T} + {T \over 2\pi} [-1 -\log{2\pi}] + O\Bigl( T\bigl({T
\over x}\bigr)^{{1/2}-\epsilon} \Bigr) + O\Bigl({T \over (\log{T})^{M-2}} \Bigr)
\end{eqnarray*}
which gives the corollary.
\section{Conclusion}
Based on (\ref{0.1}), (\ref{0.2}), Corollary \ref{corollary1} and Corollary
\ref{corollary2}, we propose the following more precise Strong Pair Correlation
Conjecture: For any small $\epsilon > 0$ and any large $A > 1$,
\[F(x,T) = \left\{ \begin{array}{ll}
{T \over 2\pi}\log{x} + {1 \over x^2}\Bigl[{T \over 2\pi}(\log{T 
\over 2\pi})^2 - 2{T \over 2\pi}\log{T \over 2\pi} \Bigr] \\
+ O(x) + O({T \over x^{{1/2}-\epsilon}}), & \raisebox{1ex}{if $1 \leq x \leq T$,
}\\
\raisebox{-1ex}{${T \over 2 \pi}\log{T \over 2\pi} - {T \over 2\pi} + O(T({T
\over x})^{{1/2}-\epsilon})$,} & \raisebox{-1ex}{if $T \leq x \leq T^{1+
\epsilon}$,}\\
\raisebox{-1ex}{${T \over 2 \pi}\log{T \over 2\pi} - {T \over 2\pi} + O(T^{1-
\epsilon_1})$,} &
\raisebox{-1ex}{if $T^{1+\epsilon} \leq x \leq T^A$.}
\end{array} \right. \]
where $\epsilon_1 > 0$ may depend on $\epsilon$, and the implicit constants may
depend on $\epsilon$ and $M$.

\end{document}